\documentclass[12pt,a4paper]{amsart}
\input xy
\xyoption{all}
\xyoption{poly}
\usepackage{amsfonts,amsmath,amssymb,indentfirst,mathrsfs,amscd}
\usepackage[mathscr]{eucal}

\let\oldmarginpar\marginpar
\renewcommand\marginpar[1]{\-\oldmarginpar[\raggedleft\footnotesize #1]%
{\raggedright\footnotesize #1}}

\date{}
\author[I. Biswas]{Indranil Biswas}

\address{School of Mathematics, Tata Institute of Fundamental
Research, Homi Bhabha Road, Bombay 400005, India}

\email{indranil@math.tifr.res.in}

\author[P.B. Gothen] {Peter Gothen}
\address{Departamento de Matematica Pura,
Facultade de Ciencias, Rua do Campo Alegre 687,
4169-007 Porto Portugal}
\email{pbgothen@fc.up.pt}

\author[M. Logares]{Marina Logares}
\address{Instituto de
Ciencias Matem\'aticas CSIC-UAM-UC3M-UCM, Serrano 113
Bis 28006 Madrid Spain}

\email{marin.logares@icmat.es}

\thanks{All three authors are members of VBAC (Vector Bundles on Algebraic Curves). Second and third authors Partially supported by CRUP through Ac\c{c}\~{a}o Integrada Luso-Espanhola nº E-38/09 and by FCT (Portugal) through the projects PTDC/ MAT/099275/2008 and PTDC/MAT/098770/2008, and the Centro de Matem\'atica da Universidade do Porto.}

\title[On moduli spaces of Hitchin pairs]{On moduli spaces of Hitchin pairs}

\DeclareMathOperator{\rk}{rk \,}          
\DeclareMathOperator{\codim}{codim\,}
	
\DeclareMathOperator{\Jac}{Jac}        
\DeclareMathOperator{\Pic}{Pic}        

\DeclareMathOperator{\Gr}{Gr\,}

\DeclareMathOperator{\img}{Im \,}           
\DeclareMathOperator{\coker}{coker\,}       

\DeclareMathOperator{\ad}{ad\,}             
\DeclareMathOperator{\Hom}{Hom\,}           
\DeclareMathOperator{\End}{End\,}           
\DeclareMathOperator{\id}{Id\,}             
\DeclareMathOperator{\tr}{Tr\,}             

\DeclareMathOperator{\sing}{Sing} 

\begin{document}

\newtheorem{thm}{Theorem}[section]
\newtheorem{prop}[thm]{Proposition}
\newtheorem{lem}[thm]{Lemma}
\newtheorem{cor}[thm]{Corollary}

\newtheorem{defn}[thm]{Definition}
\newtheorem{ex}[thm]{Example}
\newtheorem{ass}{Assumption} 

\newtheorem{rmk}[thm]{Remark}

\newcommand{\iacute}{\'{\i}} 
\newcommand{\norm}[1]{\lVert#1\rVert} 

\newcommand{\cA}{\mathcal{A}} 
\newcommand{\cC}{\mathcal{C}}
\newcommand{\cD}{\mathcal{D}}
\newcommand{\cH}{\mathcal{H}} 
\newcommand{\cG}{\mathcal{G}} 
\newcommand{\cO}{\mathcal{O}} 
\newcommand{\cM}{\mathcal{M}} 
\newcommand{\cN}{\mathcal{N}} 
\newcommand{\cP}{\mathcal{P}} 
\newcommand{\cS}{\mathcal{S}} 
\newcommand{\cU}{\mathcal{U}} 
\newcommand{\cX}{\mathcal{X}}
\newcommand{\cT}{\mathcal{T}}
\newcommand{\cJ}{\mathcal{J}}

\newcommand{\mM}{\mathscr{M}} 

\newcommand{\CC}{\mathbb{C}} 
\newcommand{\QQ}{\mathbb{Q}} 
\newcommand{\HH}{\mathbb{H}} 
\newcommand{\RR}{\mathbb{R}} 
\newcommand{\ZZ}{\mathbb{Z}} 
\newcommand{\PP}{\mathbb{P}} 
\newcommand{\VV}{\mathbb{V}} 
\newcommand{\EE}{\mathbb{E}}

\renewcommand{\lg}{\mathfrak{g}} 
\newcommand{\lh}{\mathfrak{h}} 
\newcommand{\lu}{\mathfrak{u}} 
\newcommand{\la}{\mathfrak{a}} 
\newcommand{\lb}{\mathfrak{b}} 
\newcommand{\lm}{\mathfrak{m}} 
\newcommand{\lgl}{\mathfrak{gl}} 
\newcommand{\ext}{\mathrm{ext}\,} 

\begin{abstract}
Let $X$ be a compact Riemann surface $X$ of genus at--least two.
Fix a holomorphic line bundle $L$ over $X$.
Let $\mathcal M$ be the moduli space of Hitchin pairs
$(E\, ,\phi\in H^0(End_0(E)\otimes L))$ over $X$ of rank $r$ and
fixed determinant of degree $d$. The following conditions
are imposed:
\begin{itemize}
\item ${\rm deg}(L)\,\geq\,2g-2$, $r\, \geq\, 2$, and $L^{\otimes r}
\,\not=\, K_X^{\otimes r}$,

\item $(r,d)\,=\, 1$, and

\item if $g\,=\, 2$ then $r\,\ge\, 6$, and if $g\,=\, 3$
then $r\,\ge\, 4$.
\end{itemize}
We prove that that the isomorphism class of
the variety $\mathcal M$ uniquely determines the isomorphism class of
the Riemann surface $X$. Moreover, our analysis shows that $\mathcal M$ is
irreducible (this result holds without the additional hypothesis on
the rank for low genus).
\end{abstract}

\maketitle

\section{Introduction}

The classical Torelli theorem says that the isomorphism class of a
smooth complex projective curve is uniquely determined by the
isomorphism class of its polarized Jacobian (the polarization is given
by a theta divisor). This means that if $(\Jac(X)\, , \theta)\,\cong\,
(\Jac(X')\, ,\theta')$, then $X\,\cong\,X'$.  Given any moduli space
associated to a smooth projective curve, the corresponding Torelli
question asks whether the isomorphism class of the moduli space
uniquely determines the isomorphism class of the curve. The answer is
affirmative in many situations. For instance, any moduli space of
vector bundles with fixed determinant with degree and rank coprime;
this was proved by Mumford and Newstead \cite{MN} for rank two, and
extended to any rank by Narasimhan and Ramanan \cite{NR}. They show
that the second intermediate Jacobian of the moduli space is
isomorphic to the Jacobian of the curve. Since the Picard group of the
moduli space is $\mathbb Z$, the second intermediate Jacobian has a
canonical polarization. This reduces the question to the original
Torelli theorem. This result for vector bundles has been crucial for
proving Torelli theorems for moduli spaces of vector bundles with
additional structures, for example, a Higgs field or a section (see
\cite{BiG}, \cite{M}).

Let $X$ be a compact connected Riemann surface of genus $g$, with
$g\,\geq\, 2$. Fix a holomorphic line bundle $L$ over $X$. An
\textit{$L$--twisted Higgs bundle} or a \textit{Hitchin pair}
consists of a holomorphic
vector bundle $E\,\longrightarrow\,X$ and a section $\phi\,\in \,
H^0(X, End(E)\otimes L)$. There is an appropriate notion of
(semi)stability of these objects. We recall that if $L\,=\, K_X$, a
Hitchin pair is a Higgs bundle. Higgs bundles were introduced by
Hitchin in \cite{Hi}, where the moduli spaces of Higgs bundles
were constructed using gauge theoretic methods. Nitsure in
\cite{Ni} gave a GIT
construction of a coarse moduli scheme $\cM(r,d,L)$ of
$S$--equivalence classes of semistable
Hitchin pairs of rank $r$ and degree
$d$. This moduli space is a normal quasi--projective variety.

The determinant map on $\cM(r,d,L)$ is defined as follows:
\begin{eqnarray}
\det:\cM(r,d,L)&\longrightarrow & \cM(1,d,L)\cong \Jac^{d}(X)\times H^0
(X,L)\\
\notag (E,\phi)&\longmapsto& (\Lambda^r E, \tr(\phi))\, .
\end{eqnarray}
Fix a line bundle $\xi$ over $X$ of degree $d$.
The moduli space of Hitchin pairs with fixed
determinant $\xi$ is then the preimage
$\cM_{\xi}(r,d,L)\,:=\,\det^{-1}(\xi, 0)$.

Our goal will be to prove the following.

\begin{thm}\label{thm:main}
Let $L\longrightarrow X$ be a fixed line bundle and let
$\cM_{\xi,X}(r,d,L)$ be the moduli space of Hitchin pairs of fixed
determinant $\xi\longrightarrow X$ of degree $d$. Assume that
\begin{itemize}
\item ${\rm deg}(L)\,\geq\,2g-2$, $r\, \geq\, 2$, and $L^{\otimes r}
\,\not=\, K_X^{\otimes r}$,

\item $(r,d)\,=\, 1$, and

\item if $g\,=\, 2$ then $r\,\ge\, 6$, and if $g\,=\, 3$
then $r\,\ge\, 4$.
\end{itemize}
If $\cM_{\xi,X}(r,d,L)\,\cong\, \cM_{\xi',X'}(r,d,L')$, where the
Riemann surface $X'$ and the line bundles $L'$ and
$\xi'$ are also of the above type, then $X$ is
isomorphic to $X'$.
\end{thm}

Let $\cN_{\xi}(r,d)$ (respectively, $\cN_{\xi'}(r,d)$) be the moduli
space
of stable vector bundles $E$ over $X$ (respectively, $X'$) of rank $r$
and $\Lambda^r E\, =\, \xi$ (respectively, $\Lambda^r E\, =\,
\xi'$). Theorem \ref{thm:main} is proved by showing that if
$\cM_{\xi,X}(r,d,L)\,\cong\,
\cM_{\xi',X'}(r,d,L)$, then $\cN_{\xi}(r,d)\,\cong\,\cN_{\xi'}(r,d)$.
In fact, the strategy will be to show that the open subset
$\cU\, \subset\, \cM_{\xi}(r,d,L)$ consisting of
pairs $(E\, , \phi)$ such that $E$ is
stable is actually a vector bundle over the moduli space
$\cN_{\xi}(r,d)$. This
open subset is proven to be of codimension large enough to induce
an isomorphism of the second intermediate Jacobians
$\Jac^{2}(\cM_{\xi}(r,d,L))\stackrel{\sim}{\longrightarrow}
\Jac^{2}(\cU)$. On the other hand,
$\Jac^{2}(\cU)\, \cong\, \Jac^{2}(\cN_{\xi}(r,d))$.
Then using the fact that $\text{Pic}(\cU)\,=\,\mathbb Z$
we construct a natural polarization on $\Jac^{2}(\cU)$; this is
done following the method in \cite[Section 6]{M}.
This polarization is taken to the natural polarization on
$\Jac^{2}(\cN_{\xi}(r,d))$.
This proves Theorem \ref{thm:main} using
the earlier mentioned result of \cite{MN}, \cite{NR}. The details are
given in Section~\ref{sec:torelli-theorem} below.

As a byproduct of our computations we also obtain the following
theorem (proved in Section~\ref{sec:codimc} below):

\begin{thm}\label{thm:main2}
Let $L\longrightarrow X$ be a fixed line bundle and let
$\cM_{\xi,X}(r,d,L)$ be the moduli space of Hitchin pairs of fixed
determinant $\xi\longrightarrow X$ of degree $d$. Assume that
\begin{itemize}
\item ${\rm deg}(L)\,\geq\,2g-2$, $r\, \geq\, 2$, and $L^{\otimes r}
\,\not=\, K_X^{\otimes r}$, and

\item $(r,d)\,=\, 1$.
\end{itemize}
Then the moduli space $\cM_\xi (r,d,L)$ is irreducible.
\end{thm}

Theorem \ref{thm:main2} was proved in \cite{Ni} under the assumption
that $r=2$.

\section{Hitchin pairs}
\label{sec:hitchin-pairs}

As before, $X$ is a compact connected Riemann surface of genus $g$,
with $g\, \geq\, 2$, and $ L\,\longrightarrow\, X $ is a holomorphic
line bundle. Fix integers $r\, \geq\, 2$ and $d$. We
consider Hitchin pairs $(E,\phi \in H^0(End(E)\otimes L))$ with
$\rk(E)=r$ and $\deg(E) = d$, as described in the Introduction. Recall
that the slope of $E$ is $\mu(E) = \deg(E) / \rk(E)$.

The following result is stated without proof in Remark 1.2.2 of \cite{bo}.

\begin{prop}\label{prop:univ}
  There is a universal vector bundle $\mathbb{E}$ on $\cM(r,d,L)
  \times X$ (respectively, $\cM_{\xi}(r,d,L) \times X$) whenever
  $(r,d)=1$.
\end{prop}

\begin{proof}
  By \cite{Ni}, the moduli space
$\cM(r,d,L)$ is the GIT quotient of an appropriate
  Quot--scheme $Q$ by $\mathrm{GL}_N({\mathbb C})$ (for suitable $N$),
and there is a
  universal vector bundle $\mathbb{E}$ over $X\times Q$
  (\cite[Proposition~3.6]{Ni}).  Moreover, the isotropy for the action
  of $\mathrm{GL}_N({\mathbb C})$ on a stable point of $Q$ is $\CC^*$
(the center
  of $\mathrm{GL}_N({\mathbb C})$), and $\CC^*$ is contained in the
isotropy
  subgroup of each point of $Q$. Also, the universal vector bundle
  $\mathbb{E}$ on $X\times Q$ is equipped with a lift of the action of
  $\mathrm{GL}_N({\mathbb C})$.

  There is a fixed integer $\delta$ such that for any $c$ in the
  center of $\mathrm{GL}_N({\mathbb C})$, the action of $c$ on a fiber
of
  $\mathbb{E}$ is multiplication by $c^\delta$. (As pointed out in
  Remark~1.2.3 of \cite{bo}, the fact that $\delta$ may be non-zero is
  the reason why $\EE$ does not, in general, descend to the quotient.)

  Fix a point $x_0$ of $X$.  We have two line bundles on $Q$.  The
  first line bundle $L_1$ is the top exterior power of the restriction
  of $\mathbb{E}$ to $x_0\times Q$. The second line bundle $L_2$ is
  the Quillen determinant line bundle for the family, i.e.,
  $$
    L_2 = \mathrm{Det}(\EE) = \det(R^0f_* \EE)\otimes \det (R^1f_*\EE) ^*
  $$
  where $f$ is the projection of $X\times Q$ to $Q$.

  Both these line bundles are equipped with a lift of the action of
  $\mathrm{GL}_N({\mathbb C})$: for $L_1$, any $c\in \CC^*$ acts as
multiplication by $c^r$, and
for $L_2$ any $c \in \CC^*$ acts as
multiplication by $c^e$, where
\begin{displaymath}
  e = d +r(1-g)
\end{displaymath}
is the Euler characteristic of $\EE$ restricted to a fiber of $f$.
Since $r$ and $e$ are coprime, we can express $-\delta$ (defined above) as
\begin{displaymath}
 -\delta = ar + be\, ,
\end{displaymath}
where $a$ and $b$ are integers.

Now replace the universal bundle $\mathbb{E}$ by
\begin{displaymath}
\EE' := \EE\otimes (L_1)^{\otimes a}\otimes (L_2)^{\otimes b}\, .
\end{displaymath}
It follows from our construction that $\CC^*$ acts trivially on the
fibers of $\EE'$.  Hence $\EE'$ descends to $X\times \cM(r,d,L)$.
\end{proof}

\section{Deformation theory for Hitchin pairs}

The infinitesimal deformations of a Hitchin pair $(E, \phi)$
are given by the first hypercohomology of the complex
\begin{equation}\label{defcomplex}
  C^{\bullet} = C^\bullet(E,\phi)\colon
  End(E)\,\xrightarrow{\ad(\phi)}\,
End(E)\otimes L\,
\end{equation}
where $\ad(\phi)(s)\,= \phi\circ s- (s\otimes 1_L ) \circ \phi$ (see,
for example, \cite{BiR}, \cite{bo}). Therefore,
Proposition~\ref{prop:univ} implies the following result:

\begin{prop}
  \label{prop:zariski-tangent}
  For any stable Hitchin pair $(E,\phi)\,\in\, \cM(r,d,L)$, the
  Zariski tangent space $T_{(E,\phi)}\cM(r,d,L)$ to $\cM(r,d,L)$ at
  the point $(E,\phi)$ is canonically isomorphic to
  $\HH^{1}(C^{\bullet})$.
\end{prop}

We have a long exact sequence {\setlength\arraycolsep{2pt}
\begin{eqnarray}\label{eq:long}
&0&\longrightarrow \HH^{0}(C^{\bullet})\longrightarrow
H^{0}(X,End(E))\xrightarrow{\ad(\phi)}
H^{0}(X,End(E)\otimes L)\\
\notag &&\longrightarrow \HH^{1}(C^{\bullet}) \longrightarrow
H^{1}(X,End(E))\xrightarrow{\ad(\phi)}
H^{1}(X,End(E)\otimes L)\\
\notag &&\longrightarrow \HH^{2}(C^{\bullet}) \longrightarrow 0
\end{eqnarray}
\cite[Remark 2.7]{BiR}. In particular, $\HH^{0}(C^{\bullet})$ can be
naturally identified with the space of global endomorphisms $\End(E,\phi)$.

Similarly, in the case of a Hitchin pair $(E, \phi)$ with
fixed determinant $\Lambda^r E = \xi$, the deformation complex is
\begin{equation}\label{cb0}
C^{\bullet}_0 \,:\,  End_0(E)\,
\stackrel{\ad(\phi)}{\longrightarrow}\,End_0(E)\otimes  L\, ,
\end{equation}
where $End_0(E)\, \subset\, End(E)$ is
the subbundle of rank $r^2-1$ given by the sheaf of trace--free
endomorphisms.

We shall need the following standard lemma.
\begin{lem}
  \label{lem:stable-hom}
  Let $(E,\phi)$ and $(E',\phi')$ be semistable Hitchin pairs. If
  $\Hom((E,\phi),(E',\phi')) \neq 0$, then $\mu(E) \leq \mu(E')$.

  If, moreover, $(E,\phi)$ and $(E',\phi')$ are stable, then any
  non-zero $\psi\in\Hom((E,\phi),(E',\phi'))$ is an isomorphism.
\end{lem}

\begin{proof}
  Take any non-zero $\psi\in \Hom((E,\phi),(E',\phi'))$, so
$\psi\colon E \to E'$ is a
  homomorphism such that
$$
(\psi\otimes\text{Id}_L)\circ \phi =
  \phi'\circ \psi\, \in\, H^0(X,\, Hom(E\, ,E'\otimes L))\, .
$$
Then the subsheaf $\ker(\psi) \subset E$ is
  $\phi$--invariant and the subsheaf $\mathrm{im}(\psi) \subset E'$ is
  $\phi'$--invariant. Now the lemma follows form the conditions of
semistability and stability.
\end{proof}

\begin{prop}
  \label{prop:HH}
  Let $(E,\phi)$ be a stable Hitchin pair. Then
  \begin{align*}
      \HH^0(C^{\bullet}(E,\phi)) &\cong \CC\, . \\
  \intertext{If, moreover, we assume that $L = M\otimes K$ for
  a line bundle $M$ satisfying $\deg(M) \geq 0$, then}
    \HH^2(C^{\bullet}(E,\phi)) &\cong
    \begin{cases}
      0 &\text{if $E \not\cong E\otimes M$}\, ,\\
      \CC &\text{if $E \cong E\otimes M$}\, .
    \end{cases}
  \end{align*}
\end{prop}

\begin{proof}
  Since $\HH^0(C^\bullet(E,\phi)) \cong \End(E,\phi)$, the first
  statement is immediate from Lemma~\ref{lem:stable-hom}.

  To prove the second statement, consider the by Serre duality
$\HH^{2}(C^{\bullet}(E,\phi))\,=\, \HH^0((C^\bullet)^{\vee}\otimes
K_X)^*$, where
  \begin{displaymath}
    (C^\bullet)^{\vee}\otimes K_X\colon
    End(E)\otimes M^{-1} \xrightarrow{-\ad(\phi)} End(E)\otimes
M^{-1}L\, .
  \end{displaymath}
  But $\HH^0((C^\bullet)^{\vee}\otimes K_X)$ is isomorphic to the
  space of global homomorphisms of Hitchin pairs
  $\Hom((E,\phi),(E\otimes M^{-1},\phi \otimes 1_{M^{-1}}))$ (cf.\
  \cite{gk}). Hence the second statement also follows from
  Lemma~\ref{lem:stable-hom}.
\end{proof}

For the remainder of the paper we assume that
\begin{itemize}
  \item $\deg(L) \geq 2g-2$,
  \item $L^{\otimes r}\,\not=\, K_X^{\otimes r}$, and
  \item $r$ is coprime to $d$.
\end{itemize}

\begin{prop}
  \label{prop:smooth}
All the irreducible components of the moduli spaces $\cM(r,d,L)$ and
  $\cM_\xi(r,d,L)$ are smooth.
\end{prop}

\begin{proof}
  By Lemma~\ref{lem:stable-hom}, the automorphism group of a stable
  Hitchin pair coincides with the center $\CC^*$ of
  $\mathrm{GL}_r(\CC)$. Moreover, under our assumptions,
  Proposition~\ref{prop:HH} gives the vanishing $\HH^2(C^\bullet) =
  0$. Hence the result follows from the existence of a universal
  family (Proposition~\ref{prop:univ}) and Theorem~3.1 of \cite{BiR}.
\end{proof}

\begin{prop}\label{pN}
For any stable Hitchin pair $(E,\phi)$,
$$
\dim T_{(E,\phi)}\cM(r,d,L) = r^2(\deg L)+1\, ,
$$
and every irreducible component of $\cM(r,d,L)$ is smooth of this
dimension.
\end{prop}

\begin{proof}
The Euler characteristic of the complex $C^{\bullet}$ is the following
$$
\chi(C^{\bullet})
=\dim \HH^{0}(C^{\bullet})-\dim
\HH^{1}(C^{\bullet})+\dim \HH^{2}(C^{\bullet})
=\chi(End(E))-\chi(End(E)\otimes L)\, .
$$
Hence,
\begin{equation}\label{d1}
\dim \HH^{1}(C^{\bullet})=\dim
\HH^{0}(C^{\bullet})+\dim \HH^{2}(C^{\bullet})+r^{2}\deg(L)\, .
\end{equation}
Thus, by Proposition~\ref{prop:HH} we have $\dim \HH^{1}(C^{\bullet})
= r^2(\deg L)+1$. The rest follows from
Propositions~\ref{prop:zariski-tangent} and \ref{prop:smooth}.
\end{proof}

\begin{prop}\label{prop:dimfixeddet}
  For any stable Hitchin pair $(E,\phi)$ of fixed determinant
  $\Lambda^r E = \xi$,
$$
\dim T_{(E,\phi)}\cM_\xi(r,d,L) = (r^2-1)\deg L\, ,
$$
and every irreducible component of $\cM_\xi(r,d,L)$ is smooth of this
dimension.
\end{prop}

\begin{proof}
 The proof is same as that of Proposition~\ref{pN} after considering
  the fixed determinant deformation complex (\ref{cb0}).
\end{proof}

Let
\begin{equation}\label{cU}
\cU\, \subset\, \cM_{\xi}(r,d,L)
\end{equation}
be the Zariski open subset parametrizing all pairs $(E\, ,\phi)$
such that the underlying vector bundle $E$ is stable. Let ${\mathcal
N}_\xi(r,d)$ be the moduli space of stable vector bundles $E$ of rank
$r$ with $\Lambda^r E\, =\, \xi$. We have  a forgetful map
$f:\cU\longrightarrow {\mathcal N}_\xi(r,d)$ defined by
$(E\, ,\phi)\longmapsto E$.

\begin{lem}\label{fm}
The above forgetful map
$$f:\cU\longrightarrow  {\mathcal N}_\xi(r,d)$$
makes $\cU$ an algebraic vector bundle over
${\mathcal N}_\xi(r,d)$.
\end{lem}

\begin{proof}
Take any $(E\, ,\phi)\, \in\, \cU$.
We have a map
$$
H^{0}(X, End_0(E)\otimes L)\, \longrightarrow\, f^{-1}(E)
$$
defined by $\psi\, \longmapsto\, (E\, ,\psi)$, where
$f$ is the map in the statement of the lemma.

Since $E$ is stable, $H^0(X,\text{End}(E))\,=\, {\mathbb C}\cdot
\text{Id}_E$. Hence the above map identifies the fiber $f^{-1}(E)$
with $H^{0}(X, End_0(E)\otimes L)$.

The given conditions that $E$ is stable, $\text{deg}(L)\, \geq\, 2g-2$,
and $L^{\otimes r}\, \not=\,
K^{\otimes r}_X$, imply that
$$
\dim H^{0}(End_0(E)\otimes L)\,=\, (r^2-1)\cdot(\text{deg}(L)+1-g)\, ,
$$
in particular, this dimension is independent of
$E$. Therefore, $f$ makes $\cU$ into an algebraic vector bundle over
${\mathcal N}_\xi(r,d)$.
\end{proof}

\section{Bialynicki--Birula and Bott--Morse stratifications}\label{sec:strata}





It is an important feature of the moduli space of Hitchin pairs that
it admits an action of
$\CC^{\ast}$:
\begin{eqnarray*}
\CC^{\ast}\times \cM(r,d,L)&\longrightarrow &\cM(r,d,L)\\
(t,(E,\phi))&\longmapsto &(E,t\phi).
\end{eqnarray*}
This action gives rise to a stratification of the moduli space, which
can be interpreted from a Morse theoretic point of view. Next we
recall how this comes about. We start by recalling another
interpretation of the moduli space from a gauge theoretic point of
view.

The notion of stability for a twisted Higgs bundle $(E\, ,\phi)$ is
related to the existence of a special Hermitian metric on $E$. To
explain this, fix a K\"ahler form on $X$. Let $\Lambda$ be the
contraction of differential forms on $X$ with the K\"{a}hler form. Fix
a Hermitian structure on $L$ such that the curvature of the Chern
connection (the unique connection compatible with both the holomorphic
and Hermitian structures) is a constant scalar multiple of the
K\"ahler form.  Using this Hermitian structure on $L$, the dual line
bundle $L^\vee$ is identified with the $C^\infty$ line bundle
$\overline{L}$. Let $E\, \longrightarrow\, X$ be a holomorphic
Hermitian vector bundle and $\phi$ a holomorphic section of
$\text{End}(E)\otimes L$. Using the identification of $\overline{L}$
with $L^\vee$, the adjoint $\phi^{\ast}$ is a section of
$\text{End}(E)\otimes L^\vee$.

A semistable $L$--twisted Higgs bundle is called
\textit{polystable} if it is a direct sum of
stable $L$--twisted Higgs bundles, all of the same slope. The
following Theorem is due to Li \cite{Li}; it also follows from the
general results of \cite{bradlow-garcia-prada-mundet:2003}.

\begin{thm}\label{thm:Li}
Let $(E,\phi)$ be a $L$--twisted Higgs bundle. The existence of a
Hermitian metric $h$ on $E$ satisfying
\begin{equation}\label{eq:HK}
\Lambda F_h + [\phi\, ,\phi^{\ast}]\,=\,\lambda \id\, ,
\end{equation}
for some $\lambda\,\in\, \mathbb R$, is equivalent to the polystability
of $(E,\phi)$.
\end{thm}

Here $F_{h}$ is the curvature of the Chern connection on $E$. The constant
$\lambda$ is determined by the slope of $E$. The smooth section
$[\phi\, ,\phi^{\ast}]$ of $\text{End}(E)$ is the contraction of
the section $\phi\phi^{\ast} -\phi^{\ast}\phi$ of
$\text{End}(E)\otimes L\otimes L^\vee$.

Fix a $C^\infty$ vector bundle $E\, \longrightarrow\, X$ of
degree $d$, and fix a Hermitian structure $h$ on $E$. Let $\cA$ is the
space of all unitary connections on $(E,h)$.
The equation in (\ref{eq:HK}) corresponds to the moment map for the
action of the unitary group on the product K\"{a}hler manifold
$\cA\times End(E)$ (whose K\"ahler metric is induced by the Hermitian
metrics on $E$ and $L$). The moduli space of stable $L$--twisted
Higgs bundles is then obtained as the K\"{a}hler quotient,
and hence the moduli space inherits a K\"{a}hler structure.

The restriction of the $\CC^{\ast}$-action to $S^1\subset \CC^{\ast}$
preserves the induced K\"{a}hler form on $\cM(r,d,L)$.  Thus we have a
Hamiltonian action of the circle and the associated moment map is
\begin{eqnarray*}
\mu:\cM(r,d,L)&\longrightarrow& \RR\\
(E, \phi) &\longmapsto& \frac{1}{2}||\phi||^{2}\, .
\end{eqnarray*}
It has a finite number of critical submanifolds, and $L$-twisted Higgs
bundles of the form $(V,0)$ are the absolute minima.

Let $F$ be the fixed point set for the ${\mathbb C}^*$--action on
$\cM(r,d,L)$.  This fixed point set is a disjoint union of connected
components which we denote by $F_\lambda$, so $F=\bigcup_\lambda
F_\lambda$. For any component $F_\lambda$, define
$$
U^{+}_\lambda =\{ p\in M;\, \lim_{t\to 0} tp \in F_{\lambda}\}
$$
and
$$
U^{-}_\lambda=\{ p\in M;\, \lim_{t\to \infty} tp\in F_{\lambda}\}.
$$
The sets $U^{+}_\lambda$ are strata for $\cM(r,d,L)$; the resulting
stratification is  called the
\textit{Bialynicki--Birula stratification}.

Considering the moment map $\mu$ as a Morse function, we obtain
another stratification.
To a critical submanifold $F_\lambda$ we can assign an unstable
manifold, also called \textit{upwards Morse flow},
$$
\widetilde{U}^{+}_{\lambda}=\{x\in\cM_{\xi}(r,d,L); \lim_{t\rightarrow
-\infty}\psi_{t}(x)\longmapsto F_\lambda \}\, ,
$$
where $\psi_{t}$ is the gradient flow for $\mu$.
Similarly, we have the stable manifold
$$
\widetilde{U}^{-}_{\lambda}=\{x\in\cM_{\xi}(r,d,L); \lim_{t\rightarrow
\infty}\psi_{t}(x)\longmapsto F_\lambda \}
$$
which is known as the \textit{downwards Morse flow}. Now,
$\{\widetilde{U}^{+}_{\lambda}\}_\lambda$
give a stratification, which is called the \textit{Morse
stratification}.

In \cite[Theorem 6.16]{K}, Kirwan proves that the stratifications
$\widetilde{U}^{+}$ and $U^{+}$ coincide, and similarly
$\widetilde{U}^{-}=U^{-}$ (cf.\ Hausel \cite{ha} for the first
application of this result in the context of moduli of Higgs bundles).

Let
\begin{equation}\label{e2}
h\,:\, \cM_{\xi}(r,d,L)\,\longrightarrow\,
\bigoplus_{i=2}^r H^0(X,\, L^{\otimes i})
\end{equation}
be the \textit{Hitchin map} defined by $(E\, ,\phi)\,
\longmapsto\, \sum_{i=2}^r \tr(\Lambda^i\phi)$. The inverse
image ${\mathcal H}^{-1}(0)\, \subset\, \cM_{\xi}(r,d,L)$
is called the \textit{nilpotent cone}.

The following result was observed by Hausel \cite{ha} in the context of Higgs
bundles. It generalizes to $L$-twisted Higgs bundles with essentially
the same proof.
\begin{prop}\label{prop:nilpotentcone}
The nilpotent cone coincides with the downwards Morse flow.
\end{prop}

\begin{proof}
First note that by \cite[Theorem 6.16]{K},
the downwards Morse flow coincides with
$$
\bigcup_{\lambda} \{p \in \cM_{\xi}(r,d,L);\,  \lim_{t\to \infty} tp\in
F_{\lambda} \}\, ,
$$
where $F_{\lambda}$ are the set of components of the fixed point set of
the $\CC^\ast$--action on our moduli space $\cM_{\xi}(r,d,L)$.
For the Hitchin map $h$ in \eqref{e2},
$$
h(\lim_{\lambda \to \infty} \lambda
p)\, =\, \lim_{\lambda\to\infty}(\lambda h(p))
$$
this implies that $h(p)=0$ for any point $p$ in the downwards Morse
flow. Hence, the downwards Morse flow is contained in the nilpotent
cone.

To prove the converse, recall that the
Hitchin map $h$ is proper \cite{Hi}, \cite{Ni}.
Hence  for any point $p$ of the nilpotent cone, the
$\CC^{\ast}$--orbit of $p$ is compact. So the limit $\lim_{\lambda\to
\infty}\lambda p$ exists, hence $p$ is in $F_{\lambda}$.
\end{proof}

\section{A codimension computation}\label{sec:codimc}

Let
$$
\cS=\{(E,\phi); E\; \mathrm{\, is\, not\, stable}\}\, \subset\,
\cM_{\xi}(r,d,L)
$$
be the subscheme of the moduli space $\cM_{\xi}(r,d,L)$ that
parametrizes all Hitchin pairs $(E\, ,\phi)$ such that
the underlying vector bundle $E$ is not semistable (since $r$ is coprime to
$d$, the vector bundle $E$ is semistable if and only if it is stable).
Also, consider the following set
$$
\cS'=\{(E,\phi); \lim_{t\to 0} (E,t\phi)\notin \cN_\xi(r,d)\}\, ,
$$
where $\cN_\xi(r,d)$ is the moduli space of stable vector bundles
$E\,\longrightarrow\, X$ with rank $r$ and $\Lambda^r E\cong \xi$.

\begin{prop}
  \label{prop:s-s-prime}
The equality $\cS'\,=\,\cS$ holds.
\end{prop}

\begin{proof}
Take any Hitchin pair $(E\, ,\phi)$ such that the underlying
vector bundle $E$ is stable. Then
$$
\lim_{t\to 0} (E,t\phi)\,=\, (E\, ,0)\, \in\, \cN_\xi(r,d)\, .
$$
Hence $\cS'\,\subset\, \cS$.

To prove the converse, take any $(E\, ,\phi)\,\in\, S$.
Since $E$ is not semistable, it has a unique Harder--Narasimhan
filtration
$$
E=E_m\supset E_{r-1}\supset \cdots E_{1}\supset E_0 =0\, .
$$
We recall that $E_i/E_{i-1}$ is semistable for all
$i\,\in\, [1\, ,m$, and furthermore,
$\mu(E_i/E_{i-1}) \,>\,  \mu(E_{i+1}/E_i)$. Following
Atiyah--Bott, \cite{AB}, set $D_{i}=E_{i}/E_{i-1}$,
$n_i=\rk (D_i)$, $k_i=\deg(D_i)$ and $\mu_i\,=\,
k_i/n_i$. Note that $\sum_{i=1}^m n_i= r$, and $\sum_{i=1}^m
k_i=d$. Consider the non-increasing sequence
$(\mu_1, \ldots, \mu_r)$, where each $\mu_i$ is repeated
$n_i$ times. The vector $(\mu_1, \ldots, \mu_r)$ is
called the \textit{type} of $E$.

We will analyze the limit of $(E,t\phi)$ as $t\to 0$.
For that we will recall from \cite{AB},
\cite{Sh} a partial ordering of types. For any
type $\underline{\delta}\, :=\, (\delta_1, \ldots, \delta_r)$,
consider the polygon in ${\mathbb R}^2$ traced by the points
$(\sum_{i=1}^b \delta_i\, ,b)$, $1\leq b\leq r$, and the point
$(0\, ,0)$. For another type $\underline{\delta}'$, we say
that $\underline{\delta}\,\geq\, \underline{\delta}'$ if the
polygon for $\underline{\delta}$ is above the polygon for
$\underline{\delta}'$. This ordering is complete when $r=2$.

In a family of vector bundles, the type increases under
specialization; see \cite[p.~18, Theorem 3]{Sh} for the
precise statement. Consider the map
$$
\tau\, :\, {\mathbb C}^*\, \longrightarrow\, \cM_{\xi}(r,d,L)
$$
that sends any $t$ to the point representing the Hitchin pair
$(E,t\phi)$. The Hitchin map $h$ in \eqref{e2} is proper, \cite{Hi},
\cite{Ni}, and $\lim_{t\to 0}h((E,t\phi))\,=\, 0$. Hence the above map
$\tau$ extends to a map
$$
\widehat{\tau}\, :\, {\mathbb C}\, \longrightarrow\,
\cM_{\xi}(r,d,L)\, .
$$
Let $\mathbb{E}$ be a universal vector bundle over
$\cM_{\xi}(r,d,L)\times X$ (see Proposition \ref{prop:univ}).
Consider the family of vector bundles $(\widehat{\tau}\times
\text{Id}_X)^*\mathbb{E}$ parametrized by $\mathbb C$.
For any $t\,\in\, {\mathbb C}^*$, the vector
bundle $((\widehat{\tau}\times
\text{Id}_X)^*\mathbb{E})\vert_{\{t\}\times X}\,=\, E$ is
not semistable. Hence from
\cite[p. 18, Theorem 3]{Sh} it follows that
$((\widehat{\tau}\times
\text{Id}_X)^*\mathbb{E})\vert_{\{0\}\times X}$ is not semistable.
In other words, $(E,\phi)\in \cS'$. This completes the proof.
\end{proof}

Now we are in a position to estimate the codimension of $\cS$.
We use the notation of Section \ref{sec:strata}. Let
$N\,:=\,h^{-1}(0)$ be the nilpotent cone. From
Proposition \ref{prop:nilpotentcone} we know that $N=\bigcup_\lambda F_\lambda$.

Now the complex dimension of the upwards Morse flow is:
$$
\dim (U^{+}_{\lambda})=\dim
(T\cM_{\xi}|_{F_\lambda})_{>0}+\dim(F_\lambda)
$$
because
$$
\dim (U^{+}_{\lambda})+\dim (T\cM_{\xi}|_{F_\lambda})_{<0}=\dim
\cM_{\xi}\, .
$$
Since $\cS'=\cS$, from Bott--Morse theory we know that
$\cS'=\bigcup_{\lambda\neq 0}U^{+}_{\lambda}$, so $\codim
\cS=\min_{\lambda\neq 0} \codim U^{+}_{\lambda}$. Finally
$$
\codim (U^{+}_{\lambda})=\dim (T\cM_{\xi}|_{F_{\lambda}})_{<0}
$$
which is half the Morse index\footnote{The Morse index is the real
dimension which is twice the complex dimension.}, at
$F_{\lambda}$, for the perfect Morse function $\mu$.

The critical points of the Bott--Morse function $\mu$ are exactly the
fixed points of the $\CC^{\ast}$--action on $\cM_{\xi}(r,d,L)$. If
$(E,\phi)$ corresponds to a fixed point, then it is a \emph{Hodge
  bundle}, i.e., it is of the form
\begin{equation}\label{eqn:Hodge}
E=E_{0}\oplus E_{1}\oplus \cdots\oplus E_{m},
\end{equation}
with $\phi(E_{i})\, \subset\, E_{i+1}\otimes L$ (see \cite{Hi},
\cite{simpson:1992}).  Consequently, the deformation complex in
\eqref{defcomplex} decomposes as
$$
C^{\bullet}(E,\phi)=\bigoplus_{k} C^{\bullet}_{k}(E,\phi)\, ,
$$
where for each $k$,
$$
C^{\bullet}_k (E,\phi):\bigoplus_{j-i=k}\Hom
(E_{i},E_{j})\stackrel{\ad(\phi)}{\longrightarrow}
\bigoplus_{j-i=k+1}\Hom(E_{i},E_{j})\otimes L\, .
$$
Therefore, the tangent space to $\cM_{\xi}(r,d,L)$
at the point $(E\, ,\phi)$ has a decomposition
$$
\HH^{1}(C^{\bullet}(E,\phi))= \bigoplus_{k} \HH^{1}(C^{\bullet}_{k}(E,\phi))
$$
(see Proposition \ref{prop:zariski-tangent}),
and half the Morse index at $(E,\phi)$ is
\begin{equation}\label{eq:dim}
\sum_{k>0} \dim
\HH^{1}(C^{\bullet}_{k}(E,\phi))\,
=\,\sum_{k>0}-\chi(C^{\bullet}_{k}(E,\phi))
\end{equation}
(note that the above equality follows from Proposition
\ref{prop:HH}).

We will estimate $\chi(C^{\bullet}_{k}(E,\phi))$ first,
using a similar argument to the one given in
\cite[Lemma~3.11]{bgg}. Define
$$C_k:=\bigoplus_{j-i=k}\Hom(E_{i},E_{j})$$ and
$\Phi_{k}:=\ad(\phi)|_{C_k}$. Then we have the homomorphism
$$
\Phi_{k}\, :\, C_{k}\,\longrightarrow \,C_{k+1}\otimes L\, .
$$

\begin{prop}\label{prop:eulerk}
Let $(E,\phi)$ be a stable Hitchin pair which corresponds to a fixed
point of the $\CC^{\ast}$--action. Then
$$
\chi(C^{\bullet}_{k}(E,\phi))\le (1-g) (\rk(C_k)-\rk(C_{k+1})) -
\deg(L) (\rk(C_{k+1})-\rk(\Phi_{k}))\, .
$$
\end{prop}

\begin{proof}
Note that
\begin{equation}\label{a1}
\chi(C^{\bullet}_{k}(E,\phi))=\deg(C_k)-\deg(C_{k+1})
-\rk(C_{k+1})\deg(L) + (\rk (C_{k}) -\rk (C_{k+1})) (1-g)\, ,
\end{equation}
so we will first bound $\deg(C_{k})-\deg(C_{k+1})$. For
that consider the short exact sequences
\begin{eqnarray}
0\longrightarrow \ker(\Phi_k )\longrightarrow C_k \longrightarrow
\img(\Phi_k )\longrightarrow 0\\
0\longrightarrow \img(\Phi_k )\longrightarrow C_{k+1}\otimes L
\longrightarrow \coker (\Phi_{k})\longrightarrow 0\, .
\end{eqnarray}
{}From these,
\begin{equation}\label{a2}
\deg(C_k )-\deg(C_{k+1})\,=\, \deg(\ker(\Phi_k )) +\deg(L)\rk(C_{k+1})
-\deg(\coker(\Phi_k ))\, .
\end{equation}
Clearly $\ker(\Phi_{k})\subset End_{0}(E)$. In Lemma
\ref{lem:productss} it was proved that if the pair
$(E,\phi)$ is stable, then $(End_{0}(E),\Phi_{k})$ is a
semistable pair; so we obtain
\begin{equation}\label{ineq:degle0}
\deg(\ker(\Phi_{k}))\le 0\, .
\end{equation}

In view of \eqref{a1}, \eqref{a2} and \eqref{ineq:degle0}, to
prove the proposition it suffices to show that
\begin{equation}\label{a3}
-\deg(\coker(\Phi_k))\, \geq\, -\deg(L)(\rk(C_{k+1})-\rk(\Phi_{k}))
\, .
\end{equation}

Consider the dual homomorphism $\Phi_{k}^{t}\,:\,C^{\ast}_{k+1}\otimes
L^{-1}\,\longrightarrow\, C^{\ast}_{k}$ of $\Phi_{k}$. Let
\begin{equation}\label{a5}
C_{k+1}\otimes L\,\longrightarrow\,\ker(\Phi_k^t )^{\ast}
\end{equation}
be the dual of the inclusion map $$\ker(\Phi_k^t )\,\hookrightarrow\,
C^{\ast}_{k+1}\otimes L^{-1}\, .
$$
Note that the homomorphism in \eqref{a5} vanishes identically
on the image $\text{Im}(\Phi_k)$. So we get a homomorphism
$$
f\, :\, \coker (\Phi_{k})\,\longrightarrow\, \ker(\Phi_k^t )^{\ast}
$$
which is evidently surjective.
Note that $\ker(\Phi_{k}^{t})$ is a subbundle of $C_{k+1}^{\ast}\otimes
L^{-1}$. We have a diagram
$$
\xymatrix{
0\ar[r]& \img(\Phi_k )\ar[r]\ar[rd] &C_{k+1}\otimes L\ar[r]\ar[d] &
\coker (\Phi_{k}) \ar[r]\ar[dl]^{f}& 0\\
& &(\ker(\Phi_{k}^{t}))^{\ast} & &
}
$$
and a short exact sequence
\begin{equation}\label{a6}
0\,\longrightarrow \ker(f)\,\longrightarrow
\coker(\Phi_{k})\,\longrightarrow (\ker(\Phi_{k}^{t}))^{\ast} \,
\longrightarrow 0\, .
\end{equation}
The kernel of $f$ is a torsion subsheaf of $\coker(\Phi_{k})$ (note
that $\coker(\Phi_{k})$ need not be locally free), hence from
\eqref{a6} we conclude that
\begin{equation}\label{a7}
\deg(\coker(\Phi_k))\ge \deg(\ker(\Phi_{k}^t )^{\ast})\, .
\end{equation}
As $\ker(\Phi_{k}^t )$ is a subbundle of $C^{\ast}_{k+1}\otimes
L^{-1}$, from \eqref{a7},
\begin{equation}\label{ineq:cokerker}
-\deg(\coker(\Phi_k))\le \deg(\ker(\Phi_{k}^t ))\, .
\end{equation}

We have an isomorphism $C_{k}^{\ast}\otimes L^{-1}\cong C_{-k}\otimes
L^{-1}$ and hence a commutative diagram
$$\begin{CD}
C_{k+1}^{\ast}\otimes L^{-1} @>\Phi^{t}_{k}>>  C_{k}^{\ast}\\
@V\cong VV @VV\cong V\\
C_{-k-1}\otimes L^{-1}@>-\Phi_{-k-1}\otimes 1_{L^{-1}}>>C_{-k}
\end{CD}$$
so $\ker(\Phi_{k}^{t})=\ker(\Phi_{-k-1})\otimes L^{-1}$. Hence
$$\deg(\ker(\Phi_{k}^{t}))\,
=\, \deg(\ker(\Phi_{-k-1}))-\deg(L)\rk(\ker(\Phi_{-k-1}))\, ,
$$
and then
\begin{equation}\label{ineq:kerdeg}
\deg(\ker(\Phi_{k}^{t})\,\le\, -\deg(L)\rk(\ker(\Phi_{-k-1}))\, .
\end{equation}
Notice that $\rk(\Phi_{-k-1})=\rk(\Phi_{k}^{t})=\rk(\Phi_k )$ and $\rk
(C_{k+1}) =\rk (C_{-k-1}^{\ast} )=\rk(C_{-k-1})$. Then
\begin{equation}\label{a8}
\rk(\ker(\Phi_{-k-1}))=\rk(C_{k+1})-\rk(\Phi_k )\, .
\end{equation}
Note that \eqref{ineq:cokerker}, \eqref{ineq:kerdeg} and
\eqref{a8} together imply \eqref{a3}. This completes the
proof of the proposition.
\end{proof}

The correspondence in Theorem \ref{thm:Li} gives a differential
geometric proof of the following lemma which was used above in the proof
of Proposition \ref{prop:eulerk}.

\begin{lem}\label{lem:productss}
Let $(E,\phi)$ be a stable pair, then the pair $(End(E),\ad(\phi))$ is semistable.
\end{lem}
\begin{proof}
An irreducible solution of the equations in (\ref{eq:HK}) provides us
with a semistable Hitchin pair. Now it is easy to see that the tensor
product of two solutions is again a solution to the equations, and the
same is with the dual. So if $h$ is a Hermitian structure on
$E$ that satisfies the equation in Theorem \ref{thm:Li}, then
$(End(E)=E^{\ast}\otimes E,
\ad(\phi))$ with the induced metric also satisfies the equation in
Theorem \ref{thm:Li}.
\end{proof}

\begin{prop}\label{prop:codim}
  The codimension of the subset $\cS$ of $\cM(r,d,L)$ consisting of
  pairs for which the underlying bundle is not stable, is greater than or
  equal to $(g-1)(r-1)$. Hence, it is greater than or equal to $4$
  whenever one of the following occur
\begin{itemize}
\item $g= 2$ and $r\ge 6$,
\item $g= 3$ and $r\ge 4$,
\item $ g\ge 4$ and $r\ge 2$.
\end{itemize}
\end{prop}

\begin{proof}
{}From Proposition~\ref{prop:eulerk} we need to compute the following sum
$$
\codim(\cS)
\,=\, \sum_{k>0} -\chi(C^{\bullet}_{k})
$$
$$
\ge
\sum_{k>0}\deg(L)(\rk(C_{k+1})-\rk(\Phi_{k}))+ (g-1)(\rk
(C_{k})-\rk(C_{k+1}))
$$
$$
=\, \sum_{k>0} (2g-2+l)(\rk(C_{k+1})-\rk(\Phi_{k}))+ (g-1)(\rk
(C_{k})-\rk(C_{k+1}))
$$
with $l\ge 0$ since we are assuming that $\text{deg}(L)\ge K$.
Hence,
\begin{eqnarray*}
&=\sum_{k>0}(g-1)(\rk(C_{k+1})+\rk(C_{k})-2\rk(\Phi_{k}))+ l(\rk(C_{k+1})-\rk(\Phi_{k}))\\
&\ge (g-1) \rk(C_{1})
\end{eqnarray*}
(we use that $\rk(C_{k+1})-\rk(\Phi_{k})\ge 0$).
We now just need to estimate $\rk(C_1)$, so let $r_i$ be the rank of
$E_{i}$
in the Hodge decomposition. Since $k=1$, we need to estimate $r_1 r_2
+r_2 r_3 +\cdots +r_{m-1}r_{m}$, where $m$ was the top index in
(\ref{eqn:Hodge}), but it is clearly bigger than or equal to $r-1$.
\end{proof}


\begin{proof}[Proof of Theorem~\ref{thm:main2}]
The subset $\cU$ of $L$-twisted Higgs bundles for which the
underlying bundle is stable is the complement of $\cS$. It follows
from Lemma~\ref{fm} that $\cU$ retracts onto
$\cN_\xi (r,d)$, which is well known to be connected.
Hence Proposition \ref{prop:codim} and
Proposition \ref{pN} together imply that the moduli space
$\cM_\xi (r,d,L)$ is irreducible.
\end{proof}

\section{Torelli theorem}
\label{sec:torelli-theorem}

We recall the notion of $s$--th intermediate Jacobian which will be central in our study.
Let $M$ be a complex projective manifold. Let
$$
H^{n}(M,\CC)\,=\,\bigoplus_{p+q=n}H^{p,q}(M)
$$
be the Hodge decomposition. Set $n=2s-1$, and define
$$
V_{s}\,:=\, H^{s-1,s}(M)\oplus \cdots \oplus H^{0,2s-1}
$$
Note that it can be defined for any $1\le s\le m=\dim_{\mathbb C} M$.
Then
$$
H^{2s-1}(M,\CC)\,=\,V_{s}\oplus \overline{V_{s}}\, ,
$$
and we have the projection of $H^{2s-1}(M,\CC)$ to the factor $V_s$.
Denote by $\Lambda_{s}$ the image of the composition
$$H^{2s-1}(M,\ZZ)\,\longrightarrow\, H^{2s-1}(M,\CC)
\longrightarrow V_{s}\, .
$$
The $s$--th intermediate Jacobian of $M$ is defined
to be the complex torus
$$
\Jac^{s}(M)\,:=\, \frac{V_s}{\Lambda_{s}}\, .
$$

When $s=1$, we get the Picard variety
$$
\Jac^{1}(M)=\frac{H^{0,1}(M)}{H^{1}(M,\ZZ)}=\Pic^0(M)\, .
$$
Here we will be interested in the second intermediate Jacobian
$$
\Jac^{2}(M)\,=\,\frac{H^{1,2}(M)\oplus H^{0,3}(M)}{H^{3}(M,\ZZ)}\, .
$$

A complex torus $T\,=\,V / \Lambda$, where $\Lambda$ is a cocompact
lattice in a vector space $V$, is called an \textit{abelian variety}
if it is a projective algebraic variety. The Kodaira embedding theorem
says that $T$ admits an embedding in a projective space if an only if
there exists a Hodge form on $T$, meaning a closed positive form
$\omega$ of type $(1,1)$ representing a rational cohomology class. The
conditions which determine if such a form exists are the
Riemann bilinear relations which can be formulated as follows: a class
in $H^{2}(T, \ZZ)$ is given by a bilinear form
$$
Q\,:\,\Lambda\otimes_{\ZZ}\Lambda \,\longrightarrow\, \ZZ\, , \qquad
Q(\lambda,
\lambda')\, =\, -Q(\lambda',\lambda)\, .
$$
Identifying $\Lambda_{s}\otimes_{\ZZ} \CC$ with $V_{s}\oplus
\overline{V_{s}}$, the Riemann bilinear relations are
\begin{eqnarray*}
Q(v,v')=0, \quad v,v'\in V\, ,\\
-\sqrt{-1} Q(v,v')>0, 0\neq v \in V\, .
\end{eqnarray*}
Thus, for instance, given a Hodge form $\omega$ on a complex
projective manifold $M$, the Riemann bilinear relations
$$
Q\,:\,\Lambda_{1}\otimes \Lambda_{1}\,\longrightarrow\, \ZZ
$$
given by
$$
Q(\lambda,\lambda')=\int_{M}\omega^{n-1}\wedge \lambda\wedge\lambda'
$$
produce a polarization on $\Jac^{1}(M)$ (see \cite{gh}).

We recall the basics of a mixed Hodge structure.
Let be $H$ a finite dimensional vector space over $\QQ$. A pure Hodge structure of weight $k$ on $H$ is a decomposition
$$
H_{\CC}=H\otimes \CC=\bigoplus_{p+q=k}H^{p,q}
$$
such that $H^{q,p}=\overline{H}^{p,q}$; the bar denotes complex
conjugation in $H_{\CC}$. It has two associated filtrations,
the (non-increasing) Hodge filtration $F$ on $H_{\CC}$
$$
F^{p}:=\bigoplus_{p'\ge p}H^{p',q}\, \subset\, H_{\CC}\, ,
$$
and the (non-decreasing) weight filtration $W$ defined over $\QQ$
$$
W_{m}:=\bigoplus_{p+q\le m}H^{p,q}\, .
$$

A mixed Hodge structure on $H$ consist of two filtrations: a non decreasing (weight) filtration $W$ defined over $\QQ$, and a non increasing (Hodge) filtration $F$ so that $F$ induces a Hodge filtration of weight $r$ on each rational vector space $\Gr ^{W}_{r}=W_{r}/W_{r-1}$.

Let $M_0$ be a smooth complex quasi--projective variety. The
cohomology of $M_0$ has a mixed Hodge structure \cite{De1,De2,De3}. The
above construction of second intermediate Jacobian can be
generalized to $M_0$ as follows:
\begin{equation}\label{J}
\text{Jac}^2(M_0)\, :=\, H^3(M_0,\, {\mathbb C})/
(F^2H^3(M_0,\, {\mathbb C})+ H^3(M_0,\, {\mathbb Z}))
\end{equation}
(see \cite[p. 110]{Ca}).
This intermediate Jacobian is a generalized torus \cite[p. 111]{Ca}.

Since the projection $f$ in Lemma \ref{fm} makes
$\mathcal U$ a vector bundle over ${\mathcal N}_\xi(r,d)$,
we conclude that the corresponding homomorphism
$$
f^*\, :\, H^j({\mathcal N}_\xi(r,d),\, \ZZ)
\,\longrightarrow\, H^j({\mathcal U},\, \ZZ)
$$
is an isomorphism for all $j$. In particular, it is an isomorphism
for $j\,=\, 3$. Therefore, the following proposition holds.

\begin{prop}\label{prop:bj}
$$
\Jac^{2}(\cU)\cong\Jac^{2}(\cN_{\xi}(r,d))\, .
$$
\end{prop}

\begin{lem}\label{lem:codimjac}
Let $M$ be a smooth variety and $S$ a closed subscheme of it of
codimension $k$; denote $U=M\setminus S$. Then the inclusion map
$U\,\hookrightarrow\, M$ induces isomorphism
$$
H^{j}(M,\ZZ)\cong H^{j}(U,\ZZ)
$$
for all $j<2k-1$.
\end{lem}

\begin{proof}
This lemma is proved in \cite[Lemma 6.1.1]{as}.
\end{proof}

\begin{proof}[Proof of Theorem \ref{thm:main}]
  Proposition \ref{prop:codim} and Lemma \ref{lem:codimjac} imply that
$$
\Jac^2(\cM_{\xi}(r,d,L)\,\cong\, \Jac^{2}(\cU)\, ,
$$
where $\cU$ is the open subset in Proposition \ref{prop:bj}.  Hence
from Proposition \ref{prop:bj},
$$
\Jac^2(\cM_{\xi}(r,d,L)\,\cong\, \Jac^2(\cN_{\xi}(r,d))\, .
$$
But $\Jac^2(\cN_{\xi}(r,d))\,=\,\text{Pic}^0(X)$ \cite[p. 1201,
Theorem]{MN}, \cite[p. 392, Theorem 3]{NR}.  Hence it only remains to
provide $\cM_{\xi}(r,d,L)$ with a canonical polarization, which is
done in Proposition \ref{prop:pol}. The usual Torelli theorem then
completes the proof.
\end{proof}

\section{Reconstructing the polarization}

In this section we construct a canonical polarization on
$J^{2}(\cM_{\xi}(r,d,L))$ following \cite[Section 6]{M} (see also
\cite[Section 8]{as}). To be precise, one shows that
$H^{3}(\cM_{\xi}(r,d,L))$ has a polarization which is natural in the
sense that it can be constructed when $X$ varies in a family.

\begin{prop}\label{prop:pol}
Assume that if $g=2$ then $r\ge 6$, and if $g=3$ then $r\ge 4$.
The Hodge structure $H^{3}(\cM_{\xi}(r,d,L))$ is naturally polarized. Furthermore the isomorphism
$$\Jac^2(\cM_{\xi}(r,d,L))\,\cong\, \Jac^2(\cN_{\xi}(r,d))$$
respects the polarizations.
\end{prop}

\begin{proof}
For brevity write $M=\cM_{\xi}(r,d,L)$.  We note that
Proposition~\ref{prop:codim}, Proposition~\ref{prop:bj} and
Lemma~\ref{lem:codimjac} prove that $\Pic(M)=\ZZ$, thus there is a
unique generator of the Picard group.

We make use of $\overline{M}$, the compactified moduli space of
Hitchin pairs considered in \cite{ha, Sch, Si}.  Let
$\sing(\overline{M})$ be the singular locus of the compactified moduli
space. Take a hyperplane intersection $Z$ of $M$ of codimension
$3$. By Proposition \ref{prop:codimsing} below we have $\codim
(\sing(\overline{M}))\ge 4$, hence a generic such $Z$ is smooth.  As
in the proof of \cite[Proposition~6.1]{M}, we can define
\begin{eqnarray*}\label{map:pol}
H^{3}(M)\otimes H^{3}(M)&\longrightarrow &\ZZ,\\
\beta_1\otimes \beta_2 &\mapsto&\langle\beta_1\cup\beta_2, [Z]\rangle
\end{eqnarray*}
In view of Proposition \ref{prop:codimsing}, the argument in \cite{M}
now goes through with the obvious adaptations to prove that this is a
polarization.
\end{proof}

\begin{prop}\label{prop:codimsing}
The singular locus $\sing(\overline{M})$ has codimension greater than
or equal to $(g-1)(r-1).$
\end{prop}

\begin{proof}
  Write $\overline{M} = M \cup Z$ for the compactified moduli space.
  The singular locus of $\overline{M}$ sits inside $Z$, which is an
  orbifold submanifold of codimension $1$ (Theorem 3.4 \cite{ha}). We
  can describe the singular points in $\overline{M}$ in the following way.
  Recall from \cite{ha,Sch} that $\overline{M}= (M\times \CC - N\times
  \{0\})/ \CC^{\ast}$, being $N$ the nilpotent cone defined in Section
  \ref{sec:strata}.  Singular points correspond then to fixed points
  of the $\CC^{\ast}$ action on $\cM_{\xi}(r,d,L)\times \CC$,
  i.e.\ Hitchin pairs $(E,\phi)$ for which there is a $p$--th root of the
  unity $\zeta$ such that $(E,\phi,0)\cong(E,\zeta_p \phi, \zeta_p
  \cdot 0)$. These pairs were identified by Simpson \cite{Si,Si2}:

  Let $f:E\longrightarrow E$ be the automorphism such that
  $f\phi=\zeta_p \phi f$. The coefficients of the characteristic
  polynomial of $f$ are holomorphic functions on $X$, hence constant,
  so the eigenvalues are constant. This gives a decomposition
  $E=\bigoplus_{\lambda} E_{\lambda}$ where $E_{\lambda}=\ker
  (f-\lambda)^{n}$. Since $f$ is an isomorphism $\lambda\neq 0$. Now
  $(f-\zeta_p \lambda)^{n}\phi=\zeta_{p}^{n} \phi (f-\lambda)^n$ so
  $\phi$ maps the eigenspace $E_\lambda$ to the eigenspace
  $E_{\zeta_p\lambda}$. We get then a chain of eigenvalues
  $\lambda,\zeta_{p}\lambda,\ldots,\zeta_{p}^{p-1}\lambda$ and,
  $\zeta_{p}^{p}\lambda$ becomes again an eigenvalue. It gives a
  decomposition of $(E,\phi)$ in a similar way to a Hodge bundle,
  except for that the indexing is by a cyclic group,
$$
\xymatrix@-1.1pc {
&E_{\lambda_{2}}\ar[r]&E_{\lambda_{3}}\ar[rd]& \\
E_{\lambda_{1}}\ar[ur]&&&E_{\lambda_{4}}\ar@/^3ex/@{.>}[dll]\\
&E_{\lambda_{p-1}}\ar[ul]&&}
$$
Such a Higgs bundle is called a \emph{cyclotomic Hodge bundle} by
Simpson and we write it as $(E=\oplus_{\lambda\in C_{p}},
\phi_{\lambda})$, where $C_{p}$ denotes the cyclic group of order
$p$.

We can consider such pairs as forming subvarieties of $M$. In order to
estimate their codimension, we study their space of deformations.
The deformation complex (\ref{defcomplex}) gives us the following
deformation complex for a cyclotomic Hodge bundle $(E=\oplus_{\lambda\in
  C_{p}}, \phi_{\lambda})$:
$$C^{\bullet}_{\lambda_{k}}:\bigoplus_{\lambda_{k}=\lambda_{j}-\lambda_{i}}\Hom(E_{\lambda_{i}},E_{\lambda_{j}})\longrightarrow \bigoplus_{\lambda_{k+1}=\lambda_{j}-\lambda_{i}}\Hom(E_{\lambda_{i}},E_{j}).$$
As in Section \ref{sec:codimc}, the tangent space at those points has
a decomposition
$$\HH^{1}(C^{\bullet}(E,\phi))=\HH^{1}(C^{\bullet}_{\lambda_{k}}(E,\phi)).$$
The dimension of the singular locus satisfies
$$\dim (\sing(\overline{M}))\le\dim T_{\sing(\overline{M})}=\dim \HH^{1}(C^{\bullet}_{e}),$$
where $e$ is the neutral element in $C_{p}$. Therefore the codimension
satisfies
$$\codim (\sing(\overline{M})\ge \dim \HH^{1}( C^{\bullet}_{\lambda\neq e}).$$
The computations in Proposition \ref{prop:eulerk} and Proposition
\ref{prop:codim} hold again for the cyclotomic Hodge bundles. This gives
the estimate for the codimension of the singular locus.
\end{proof}


\end{document}